\documentclass[a4paper,12pt]{article}
\usepackage[russian]{babel}
\usepackage{amsmath}
\usepackage{ulem}
\usepackage{calc}
\usepackage{vmargin}
\usepackage{amstext}
\usepackage{graphicx}
\usepackage{amsfonts}
\usepackage{amssymb}

\numberwithin{equation}{section}

\begin{document}

\begin{center}
{\bf Passive Collections of Differential Power Series Algebra}\\
O.V. Kaptsov

Institute of Computing Modeling SD RAS,\\
Academgorodok 50, 660036 Krasnoyarsk, Russia\\
E-mail: profkap@gmail.com
\end{center}












\begin{center} {\bf Abstract} \end{center}

We consider a differential algebra $\mathcal{F}$ of formal power series in infinitely many variables.
We define the important notions of a normalized set of generators for an ideal of $\mathcal{F}$ and a regular quotient algebra. The concept of the passive collection being analog of involutive system of the nonlinear partial differential equations and Gr\"obner basis  is introduced. Using special partitions of the algebra $\mathcal{F}$ and action of a semigroup on $\mathcal{F}$,
the sufficient conditions of passivity are obtained. The proof of the main theorem essentially differs from similar one for the differential equations and we obtain more general statement.


\section{Introduction}
At the end of the 19th century, S. Lie suggested to consider partial differential equations as manifolds embedded in some spaces.
Using the Galois theory and this idea, he created the theory of continuous symmetry and applied it to the research of differential equations. In the second half of the XX century, the fibre bundles and homological algebra began to be used in the geometry of differential equations \cite{Spencer,Pommaret,Seiler}. As a result, the theory became more invariant and difficult. This approach led to new problems, for example, some important differential equations are not smooth manifolds. It is therefore useful to try to apply the local methods in the study of the geometry of differential equations. One may consider algebraic and complex analytic geometries as the standard patterns \cite{Shafarevich,Abhyankar}. In algebraic geometry there is an important distinction between singular and non-singular points. Each point of affine algebraic variety corresponds to a local ring. Recall that  the local rings of non-singular points are regular.

In order to study local properties of partial differential equations we consider a differential algebra $\mathcal{F}$ of formal power series in infinitely many variables $x_1,\dots,x_n, u^i_\alpha$ ($i=1,\dots,m$, $\alpha \in \mathbb{N}^n$, where $\mathbb{N}$ are the natural numbers) with coefficients in a field $\mathbb{K}$ of characteristic $0$.
We are interested in analyzing special generators of differential ideals in $\mathcal{F}$ and the regular quotient algebras $\mathcal{F}/I$.

This paper is structured as follows. In the next Section we define the notions of normalized set of generators for an ideal $I \subset \mathcal{F}$ and the regular quotient algebra $\mathcal{F}/I$. Note that $\mathcal{F}$ is an infinite-dimensional ring and one cannot use the usual definition of the regular local ring. We prove that the remainder of the division of arbitrary series $f$ by normalized set
is uniquely determined. It follows that if an ideal $I$ has the normalized set of generators then the quotient algebra $\mathcal{F}/I$ is regular. In Section 3 we introduce the notion of a partition of a set $A$ being compatible with the action of a semigroup $G$ on $A$ and derivation operators on $\mathcal{F}$. These derivation operators induce an action of the semigroup $\mathbb{G}=\mathbb{N}^n \setminus  0$  on $\mathcal{F}$. Some compatibility criteria of  the partitions of $\mathcal{F}$ with the action $\mathbb{G}$ on $\mathcal{F}$  are given.
In the final Section we define a passive collection of the algebra $\mathcal{F}$ and determine sufficient conditions for a subset $S$ of
$\mathcal{F}$ to be the passive collection. One might think of the passive collection as an analogue of the Gr\"obner basis for an ideal of the ring of polynomials in several variables \cite{Cox}. Roughly speaking, we prove that if a partition of $\mathcal{F}$ is compatible with the action $\mathbb{G}$ on $\mathcal{F}$ and some compatibility conditions for $S\subset \mathcal{F}$ are satisfied then $S$ is passive.
It can be mentioned that our proof and classical one for partial differential equations \cite{Ritt} are essentially distinct; as well we use
slightly weaker conditions.

\section{Normalized Systems of Generators}

Let $\mathbb{K}$ be a field and let $Y$ be an arbitrary non-empty set. Consider a finite subset $W$ of $Y$ and the commutative formal power series algebra  $\mathbb{K}[[W]]$. Then
$$ \mathcal{F} =  \bigcup_{W\subseteq Y} \mathbb{K}[[W]]$$
is obviously a $\mathbb{K}$-algebra, denoted by  $\mathbb{K}||Y||$.

\noindent
{\bf Definition.} A subset $B$ of $\mathcal{F}$ is called a normalized set if there exist a unique set $\mathcal{L}\subset Y$, such that the following hold:

\noindent
1. An arbitrary series $f\in B$ is of the form $f=y+g$, where $y\in \mathcal{L}$, $g$ is a non-unit of $\mathbb{K}||\mathcal{P}||$, and
$\mathcal{P}=Y\setminus \mathcal{L}$.

\noindent
2. For every $y\in \mathcal{L}$ there is a unique $f\in B$  such that $f=y+g$, with $g\in \mathbb{K}||\mathcal{P}||$.

\noindent
The elements of $\mathcal{L}$, $\mathcal{P}$ are called leading and parametric variables respectively.

\noindent
{\bf Lemma 1.} Let $B$ be a normalized subset in $\mathcal{F}$ and let $\mathcal{P}$ be the set of corresponding parametric variables.
Then for any $f\in\mathcal{F}$, depending on $p$ ($p\geq 0$) leading variables, there exist series $f_1,\dots, f_p\in B$,
 $q_1,\dots,q_p\in \mathcal{F}$, and unique series $r\in \mathbb{K}||\mathcal{P}||$ such that

\begin{equation}\label{eq:devision}
 f = \sum_{i=1}^{p} q_i f_i + r .
\end{equation}

Proof. If $f$ is in the algebra $\mathbb{K}||\mathcal{P}||$, then $r=f$. Suppose that $f$ depends on the leading variables $y_1,\dots,y_p$ and
 $f_i = y_i + g_i \in B$, $i=1,\dots,p$. Let $f,f_1,\dots,f_p$ depend on $y_1,\dots,y_n$. Using the Weierstrass division formula \cite{Zariski} we have
$$f = q_1 f_1 +r_1(y_2,\dots,y_n) $$
with $q_1\in \mathcal{F}$. Again the Weierstrass theorem (for $r_1$) implies
$$f = q_1 f_1 +q_2 f_2+r_2(y_3,\dots,y_n) .$$
Continuing in this fashion we obtain the formula (\ref{eq:devision}),
where $r$ depends only on the parametric variables $y_{p+1},\dots,y_n$.

 To prove uniqueness, we assume that there is another representation
 $$ f = \sum_{i=1}^{p} \tilde{q_i} f_i + \tilde{r}$$
 with  $\tilde{r}$ depending on $y_{p+1},\dots,y_n$.
It follows from the last representation and (\ref{eq:devision}) that
\begin{equation}\label{eq:mu}
\sum_{i=1}^{p}\mu_i f_i = h ,
\end{equation}
where $\mu_i =(q_i - \tilde{q_i})$, $h=\tilde{r} - r$.

Consider the ring $\mathbb{K}[[y_1,\dots,y_n,y_1^{'},\dots,y_p^{'}]]$ and the equations
$$ y_1^{'} - f_1= \cdots = y_p^{'}- f_p=0.$$
The Jacobian $|\frac{\partial(f_1,\dots,f_p)}{\partial(y_1,\dots,y_p)}(0)|$ is equal to unity in $\mathbb{K}$.
 According to the implicit function theorem \cite{Abhyankar} $y_1,\dots,y_p$ are given by power series in
 $y_1^{'},\dots,y_p^{'},y_{p+1},\dots,y_n $.
Then the relation (\ref{eq:mu}) is written as
\begin{equation}\label{eq:mu2}
\sum_{i=1}^{p}\tilde\mu_i(y^{'}) y_i^{'} = h(y_{p+1},\dots,y_n) ,
\end{equation}
with $y^{'}= (y_1^{'},\dots,y_p^{'},y_{p+1},\dots,y_n)$. 
The right hand side of (\ref{eq:mu2}) is equal to zero because the left hand side is either zero or depends on at least one of
$y_1^{'},\dots,y_p^{'}$.

\noindent
{\bf Definition.} Under the conditions of Lemma 1, the series $r$  is called the remainder of the division $f$ by $B$. If $r=0$,
we say $B$ divides $f$.

\noindent
{\bf Definition.} If an ideal $I$ in $\mathcal{F}$ is generated by a normalized set $B$ then $B$ is called a normalized set of generators for $I$ .

\noindent
{\bf Definition.} Let $I$ be an ideal in  $\mathcal{F}=\mathbb{K}||Y||$. We say that the quotient algebra $\mathcal{F}/I$ is regular if it is isomorphic to $\mathbb{K}||Y^{'}||$ with $Y^{'}\subseteq Y$.

\noindent
{\bf Proposition 1.}  Let $I$ be an ideal in $\mathcal{F}$ and assume that $I$ has a normalized set of generators $B$. Then

(a) the quotient algebra $\mathcal{F}/I$ is isomorphic to $\mathbb{K}||\mathcal{P}||$, where $\mathcal{P}$ is the set of corresponding parametric variables,

(b) $\mathcal{F}$ is the direct sum of $I$ and $\mathbb{K}||\mathcal{P}||$  as vector spaces over $\mathbb{K}$.

\noindent
Proof. Let $\pi$ be the canonical projection homomorphism $\mathcal{F} \rightarrow \mathcal{F}/I$. It follows from Lemma 1 that
for every $f\in \mathcal{F}$ we have $\pi(f)=\pi(r)$, where $r$ is the remainder of the division $f$ by $B$. It is easy to see that
the map $\phi : \mathcal{F}/I \rightarrow \mathbb{K}||\mathcal{P}||$ given by the formula
$$\phi(\pi(r))=r$$
 is an isomorphism of algebras. The second assertion is an immediate consequence of Lemma 1.

\noindent
{\bf Corollary.} Let $B_1, B_2$ be two normalized sets of generators for an ideal $I\subset \mathcal{F}$ and let $\mathcal{L}_1, \mathcal{L}_2$ be the sets of corresponding leading variables. If $\mathcal{L}_1=\mathcal{L}_2$, then $B_1 =B_2$.

\noindent
Proof.  Suppose that $f_1=y+h_1\in B_1$ and $f_2=y+h_2\in B_2$. where $h_1, h_2\in \mathbb{K}[[\mathcal{P}]]$, $\mathcal{P}$ is the set of parametric variables. Hence the difference $f_1 - f_2 = h_1 - h_2$ belongs to $\mathbb{K}||\mathcal{P}||$. On the other hand, this difference lies in $I$. Since $\mathbb{K}||\mathcal{P}|| \cap I = 0$, it follows that $h_1=h_2$.

\section{Differential Power Series Algebra}

We now assume throughout that the field $\mathbb{K}$ has characteristic $0$. Denote by $\mathbb{N}$ the set of all nonnegative integers.
Then $\mathbb{N}^n$ is a monoid with generators $e_1=(1,0,\dots,0)$, \dots, $e_n=(0,\dots,0,1)$.

Henceforth $Y$ is the disjoint union of two sets
\begin{equation}\label{eq:U}
X=\{x_1,\dots,x_n\} , \qquad U= \{u^i_{\alpha} : i=1,\dots,m ; \alpha \in \mathbb{N}^n\} .
\end{equation}
Then differential operators on $\mathcal{F}=\mathbb{K}||Y||$ are given by
\begin{equation} \label{eq:D_i}
D_i f = \frac{\partial f}{\partial x_i} + \sum_{j=1,\dots,m, \ \alpha\in N^n}
\frac{\partial f}{\partial u^j_\alpha} u^j_{\alpha+e_i}
\end{equation}
with the partial derivatives $\frac{\partial f}{\partial x_i}$, $\frac{\partial f}{\partial u^j_\alpha}$
defined as usual \cite{Zariski}. The algebra $\mathcal{F}$ with the operators (\ref{eq:D_i}) is a differential power series algebra \cite{Kolchin}.
Recall some additional terminology associated to semigroups.

{\bf Definition.} A semigroup $G$ is said to act on the set $A$ if there
exists a map $(g, a) \rightarrow ga$ of $G \times A$ into $A$ satisfying
 $$g_1(g_2 a) = (g_1g_2) a , \qquad a\in A, g_1,g_2 \in G.$$

Suppose G acts on A. Then we define $T_g$ to be the map $a \rightarrow ga$,
$a \in A$. According to the preceding definition, we have
\begin{equation}\label{eq:Tg}
T_{g_1}\circ T_{g_2} = T_{g_1 g_2} .
\end{equation}
The family $\{T_g\}_{g\in G}$ is a transformation semigroup.

\noindent
{\bf Definition.} Let $\sim$ be an equivalence relation on a set $A$. We say a semigroup $G$ acting on $A$ preserves the equivalence relation, if
$$a\sim b \Rightarrow g a \sim g b, \qquad \forall g \in G.$$

It is also convenient to write the last condition as
\begin{equation}\label{eq:sim}
a\sim b \Rightarrow T_g(a) \sim T_g(b), \qquad \forall g \in G.
\end{equation}
.

Throughout this article a quotient set of $A$ relative to an equivalence relation $\sim$ is denoted by $A/\sim$ and an equivalence class of $a \in  A$ is denoted as $[a]$.

If a semigroup $G$ acting on $A$ preserves an equivalence relation $\sim$, then an action of $G$ on $A/\sim$ is given by
$$g[a] = [ga] .        $$

Henceforth we adopt the convention that $(\Gamma,\leq)$ is  a well-ordered set. Let $\{A_{\gamma}\}_{\gamma\in\Gamma}$ be  a partition of $A$ then (as is well known) there is an equivalence relation on $A$ whose equivalence classes are precisely the sets $A_{\gamma}$:
\begin{equation}\label{eq:sim2}
a\sim b \quad \Longleftrightarrow \quad \exists \gamma (a,b \in A_\gamma) .
\end{equation}
We define a well-order $\preceq$ on $A/\sim$ by
\begin{equation}\label{eq:sim2}
A_{\gamma}\preceq A_{\gamma^{'}} \quad \Longleftrightarrow \quad \gamma\leq \gamma^{'}.
\end{equation}
We will write $A_{\gamma}\prec A_{\gamma^{'}}$ if $A_{\gamma}\preceq A_{\gamma^{'}}$ and $A_{\gamma}\neq A_{\gamma^{'}}$.

\noindent
{\bf Definition.} Let $\{A_{\gamma}\}_{\gamma\in\Gamma}$ be  a partition of $A$, a semigroup $G$ acts on $A$ and preserves the equivalence relation $\sim$. We say that the partition is compatible with the action of $G$ on $A$ if $G$ acting on $A/\sim$ preserves the
well-order (\ref{eq:sim2}) and its direction, i.e.,  the following properties hold:
\begin{equation}\label{eq:order1}
(a) \qquad A_{\gamma}\prec A_{\gamma^{'}} \quad \Longrightarrow \quad [T_g (A_{\gamma})] \prec [T_g (A_{\gamma^{'}})]
\end{equation}
\begin{equation}\label{eq:order2}
(b)\qquad A_{\gamma}\prec [T_g (A_{\gamma})]  \qquad \forall \gamma\in \Gamma ,
\end{equation}
for all $g \in G$.

\noindent
{\bf Proposition 2.} Let $\{A_{\gamma}\}_{\gamma\in\Gamma}$ be  a partition of $A$ and a semigroup $G$ acts on $A$. Assume the properties (\ref{eq:sim}), (\ref{eq:order1}), (\ref{eq:order2})
 are satisfied for generators of $G$, then the partition is compatible with the action of $G$ on $A$

Proof. Every element of $G$ is a finite product of generators. Thus it suffices to prove that if $g_1$ and $g_2$ satisfy the properties of (\ref{eq:sim}), (\ref{eq:order1}), (\ref{eq:order2}) then this is also true for product $g_1g_2$. It can be deduced directly from (\ref{eq:Tg}).

Let $\{U_\gamma\}_{\gamma \in \Gamma}$ be a partition of the set $U$  (\ref{eq:U}). Consider sets
$$ Y_{\gamma} = X\cup (\bigcup_{\gamma^\prime \leq \gamma} U_{\gamma^\prime})  $$
and  the corresponding formal power series algebras
\begin{equation}\label{eq:Fgamma}
\mathcal {F}_{\gamma} = \mathbb{K}||Y_{\gamma}||.
\end{equation}
The families $\{Y_{\gamma}\}_{\gamma \in \Gamma}$ and $\{\mathcal{F}_\gamma\}_{\gamma \in \Gamma}$ are well-ordered by inclusion relation.

Let $\hat{\mathcal{F}}$ denote the difference of the sets $\mathcal{F}$ and $\mathbb{K}[[X]]$:
\begin{equation}\label{eq:F^}
\hat{\mathcal{F}} = \mathcal{F} \setminus \mathbb{K}[[X]] .
\end{equation}
The family $\{\mathcal{F}_\gamma\}_{\gamma \in \Gamma}$, where $\mathcal{F}_\gamma$ given by (\ref{eq:Fgamma}), determines a partition
$\{\Phi_\gamma\}_{\gamma \in \Gamma}$ of $\hat{\mathcal{F}}$ into the blocks
\begin{equation}
\label{eq:Phi} \Phi_\gamma = \mathcal{F_\gamma} \smallsetminus (\bigcup_{\gamma^\prime <  \gamma} \mathcal{F_{\gamma^\prime}} ) .
\end{equation}

We denote by $\mathbb{N}_{-0}^n$ the semigroup $\mathbb{N}^n\setminus 0$, where $0$ is a unit element of $\mathbb{N}^n$. Let $D^\alpha$ be the product $D_1^{\alpha_1}\cdots D_n^{\alpha_n}$ with $\alpha = (\alpha_1,\dots,\alpha_n) \in \mathbb{N}_{-0}^n$.
The semigroup $\mathbb{N}_{-0}^n$ acts on $U$ and $\hat{\mathcal{F}}$ by
\begin{equation}\label{eq:Du}
\alpha u^i_\beta = u^i_{\alpha + \beta} ,  
\end{equation}
\begin{equation}\label{eq:Daf}
\alpha f = D^\alpha (f) .
\end{equation}

\noindent
{\bf Lemma 2.} Let $\{U_\gamma\}_{\gamma \in \Gamma}$ be a partition of $U$ which is compatible with the action (\ref{eq:Du}) of the semigroup $\mathbb{N}_{-0}^n$ on $U$. Then the corresponding partition  $\{\Phi_\gamma\}_{\gamma \in \Gamma}$ of $\hat{\mathcal{F}}$ into the blocks (\ref{eq:Phi}) is compatible with the action (\ref{eq:Daf}) of $\mathbb{N}_{-0}^n$ on $\hat{\mathcal{F}}$.

\noindent
Proof. We will use Proposition 2 for generators $e_1,\dots,e_n$ of the semigroup $\mathbb{N}_{-0}^n$.
Let us first prove that $D_i$ preserves the equivalence relation on $\hat{\mathcal{F}}$, i.e., if $f_1, f_2\in \Phi_\gamma$
then there exists $\gamma^{'}\in \Gamma$ such that $D_i (f_1), D_i(f_2)\in \Phi_{\gamma^{'}}$.

Since $f_1, f_2\in \Phi_\gamma$, there exists elements $u^j_\alpha, u^l_\beta \in U_\gamma$ such that
\begin{equation}\label{eq:f1f2}
 \frac{\partial f_1}{\partial u^j_\alpha } \neq 0, \quad  \frac{\partial f_2}{\partial u^l_\beta } \neq 0 .
\end{equation}
By assumption, for all $u^j_\alpha, u^l_\beta \in U_\gamma$ there is some $\gamma^{'}$ such that
$u^j_{\alpha+e_i}, u^l_{\beta+e_i} \in U_{\gamma^{'}}.$ It follows that
$$\frac{\partial f_1}{\partial u^j_\alpha } u^j_{\alpha+e_i} \in \Phi_{\gamma^{'}},  \qquad
\frac{\partial f_2}{\partial u^l_\beta } u^l_{\beta+e_i} \in \Phi_{\gamma^{'}} .    $$
According to the relation (\ref{eq:D_i}),  $D_i f_1$ and $D_i f_2$ also lie in $\Phi_{\gamma^{'}}$.

We now must prove that $D_i$ satisfies the property (\ref{eq:order1}), i.e.,
$$ [f_1] \prec [f_2] \Rightarrow  [D_i f_1] \prec [D_i f_2]   $$
with $f_1\in \Phi_{\gamma_1}$, $f_2\in \Phi_{\gamma_2}$ ($\gamma_1 < \gamma_2$).
If $f_1\in \Phi_{\gamma_1}$ and $f_2\in \Phi_{\gamma_2}$ then there are $u^j_\alpha \in U_{\gamma_1}$ and $u^l_\beta \in U_{\gamma_2}$ such that (\ref{eq:f1f2}) holds. By hypothesis, for all $u^j_\alpha \in U_{\gamma_1}, u^l_\beta \in U_{\gamma_2}$
there exist $\gamma_3, \gamma_4 \in \Gamma$ ($\gamma_3 < \gamma_4$ ) such that
$u^j_{\alpha+e_i}\in U_{\gamma_3}, u^l_{\beta+e_i} \in U_{\gamma_4}.$ Thus we have
$$\frac{\partial f_1}{\partial u^j_\alpha } u^j_{\alpha+e_i} \in \Phi_{\gamma_3},  \qquad
\frac{\partial f_2}{\partial u^l_\beta } u^l_{\beta+e_i} \in \Phi_{\gamma_4} .    $$
In this way we see that $D_i f_1\in \Phi_{\gamma_3}$ and $D_i f_2 \in \Phi_{\gamma_4}$, according to the definition above.
Similarly, it is possible to check the relation (\ref{eq:order2}).

Let $P$ be the commutative polynomial algebra $\mathbb{K}[X_1,\dots,X_n]$. We shall use the notation $\mathbb{K}\langle U\rangle$ to denote a vector space  over $\mathbb{K}$ spanned by $U$.
We can define a left $P$-module structure on the vector space $\mathbb{K}\langle U\rangle$ as follows.
The product $p = \sum_{\alpha \in \Delta} a_\alpha X^\alpha \in P$ and $u^j_\beta \in U$  is given by the formula
$$ p u^j_\beta =  \sum_{\alpha \in \Delta} a_\alpha u^j_{\alpha+\beta}, $$
where $\Delta$ is a finite subset in $\mathbb{N}^n$, and can be extended to $\mathbb{K}\langle U\rangle$ by linearity.

\noindent
{\bf Proposition 3.} The left $P$-modules $\mathbb{K}\langle U\rangle$ and $P^m$ are isomorphic.

\noindent
Proof. The left $P$-modules $\mathbb{K}\langle U\rangle$, $P^m$ have the bases  $u^1_0, \dots, u^m_0$ and $e_1,\dots, e_m$ respectively.
Then a map $\phi$, defined by $\phi(u^i_0)=e_i$ on the basis and then extended by linearity on $\mathbb{K}\langle U\rangle$, is an isomorphism.

\noindent
{\bf Definition.} Let $V$ be a set consisting of an elements $v_1,\dots,v_k\in U$. 
Every $k$-tuple $(d_1,..., d_k)$ of elements of $P$ such that
$$d_1 v_1 + \cdots + d_k v_k = 0$$
is called a syzygy of $V$.

The collection of all syzygies of a set $V$ forms a submodule of the module $P^k$ and is denoted by $Syz V$.

If $\alpha$ and $\beta$ are elements of $\mathbb{N}^n$, we can define their product:
\begin{equation}
\label{eq:diamond} \alpha \diamond \beta = (\mu_1,\dots,\mu_n), \qquad \mu_i = max(\alpha_i, \beta_i) - \alpha_i .
\end{equation}
\noindent
{\bf Example.} Assume $v_1 =u^i_\alpha$ and $v_2 = u^i_\beta$ are elements of some set $V$ consisting of $k$ elements of $U$. In this case the
$k$-tuple $d =(X^\beta,-X^\alpha,0,\dots,0)$ is a syzygy of $V$. Moreover,
$\tau_{12} = X^{\alpha \diamond \beta}e_1 - X^{\beta \diamond \alpha}e_2$
is also a syzygy of $V$.

It is easy to generalize the last example. Consider the set
\begin{equation}
\label{eq:V} V = \{ u^i_\alpha \in U : i= 1,\dots, \bar m \leq m; \alpha \in \Delta \},
\end{equation}
where $\Delta$ is a finite subset of $\mathbb{N}^n$, and a partition of $V$ into subsets
$$ V_q = \{ u^q_\alpha \in V :  1\leq q\leq \bar m \leq m\}.$$
Suppose $V$ consists of $k$ elements, one can constitute a $k$-tuple $w\in U^k$ as follows: at first we put the elements of $V_1$, next put
the elements of $V_2$, and so on.
If $u^q_\alpha$ and $u^q_\beta$ are $i^{th}$ and $j^{th}$ components of the $k$-tuple $w$, then the syzygy
\begin{equation}
\label{eq:Syzygy} \tau_{ij} = X^{\alpha \diamond \beta}e_i - X^{\beta \diamond \alpha}e_j
\end{equation}
corresponds to the elements $u^q_\alpha , u^q_\beta$.

\noindent
{\bf Proposition 4.} Let $V$ be the finite set (\ref{eq:V}). Then the module $Syz V$ is generated by the syzygies  (\ref{eq:Syzygy}).

\noindent
The proof follows at once from the corresponding property of syzygies of monomial submodules \cite{Eisenbud} of the $P$-module $P^k$
 and Proposition 3.

\section{Passive Collections}

M\'{e}ray and Riquier were possibly first who started studying passive systems of partial differential equations \cite{Meray}.
The summary of the classical theory of these systems can be found in \cite{Ritt}.
In this section we shall introduce the notion of passive collections of the algebra $\mathcal{F}$ and
determine sufficient conditions for a subset to be a passive collection.

\noindent
{\bf Definition.} The set
$$ O(S) = \{D^\alpha s:  \alpha\in \mathbb{N}^n; s\in S\} $$
is called the orbit of $S\subset \mathcal{F}$ under $\mathbb{N}^n$.

Assume $f\in \mathcal{F}$ has the form
\begin{equation} \label{eq:f=u^i_alpha}
 f = u^i_\alpha + h^i_\alpha,
\end{equation}
with $u^i_\alpha \in U$, then the element $u^i_\alpha$ is denoted by $st f$.
If a subset $S$ of $\mathcal{F}$ consists of power series (\ref{eq:f=u^i_alpha}), then $st S= \{st f: f\in S\}.$

\noindent
{\bf Definition.} A collection $S=\{f_1,\dots,f_k\} \subset \mathcal{F}$ is said to be passive if it satisfies the following conditions:

(i) every $f \in S$ is of the form (\ref{eq:f=u^i_alpha}),
where $h^i_\alpha$ is not dependent on $u^i_\alpha$.

(ii) The ideal $\langle O(S) \rangle$ has a normalized set of generators $B$ such that
\begin{equation} \label{eq:L=O}
\mathcal{L} = O(st S) ,
\end{equation}
where $\mathcal{L}$ is a set of leading variables for $B$.

We further assume that some family $\{\mathcal{F}_\gamma\}_{\gamma \in \Gamma}$ of algebras (\ref{eq:Fgamma}) and the corresponding partition $\{\Phi_\gamma\}_{\gamma \in \Gamma}$ of $\hat{\mathcal{F}}$ (\ref{eq:F^}) into the blocks (\ref{eq:Phi}) are given.
Moreover, one writes $\Phi_{\gamma}\prec \Phi_{\gamma^{'}}$ whenever $\gamma \prec \gamma^{'}$.
It is obvious that $\langle O(S) \rangle$  is a differential ideal of the differential algebra $\mathcal{F}$.

\noindent
{\bf Definition.} A subset $S$ of $\mathcal{F}$ is said to be conditionally solvable if every power series $f\in S$ has the form
(\ref{eq:f=u^i_alpha}) and satisfies the property
\begin{equation} \label{eq:f=h<u^i}
[h^i_\alpha] \prec [u^i_\alpha].
\end{equation}
The element $u^i_\alpha$ is called the leading term of $f$ and is denoted by $lt f$.

 We define the element       
$$\theta = min\{\gamma \in \Gamma : O(S)\cap \mathcal{F}_\gamma \neq \emptyset \}.$$
and the sets
\begin{equation} \label{eq:Ogamma}
lt S = \{lt f : f\in S\},\qquad
O_\gamma = O(S)\cap \mathcal{F}_\gamma, \qquad  C_\gamma = O(S)\cap \Phi_\gamma ,  \qquad\forall \gamma \in \Gamma .
\end{equation}
Denote by $\langle O_\gamma \rangle_{\mathcal{F}_\gamma}$ an ideal of the algebra $\mathcal{F}_\gamma$ generated by $O_\gamma$.

\noindent
{\bf Definition.} Let $S=\{f_1,\dots,f_k\}$ be a conditionally solvable subset of $\mathcal{F}$ and let $d=(d_1,\dots,d_k)$ be a syzygy of $lt S$. We shall say $S$ satisfies $d$-compatibility condition if there exists $\gamma \in \Gamma$ such that
\begin{equation} \label{eq:compatib}
\sum\limits_{i=1}^{k} d_i f_i  \in  \Phi_\gamma \cap \langle O_\gamma \rangle_{\mathcal{F}_\gamma} .
\end{equation}

The main result of this section is the following theorem.

\noindent
{\bf Theorem.} Let $S=\{f_1,\dots,f_k\}$ be a conditionally solvable subset of $\mathcal{F}$. Let $\{\tau_{ij}\}$ be the generators (\ref{eq:Syzygy}) of the syzygy module $Syz\ lt S$.
Suppose the partition $\{\Phi_\gamma\}_{\gamma \in \Gamma}$ of $\hat{\mathcal{F}}$ is compatible with the action (\ref{eq:Daf}) of $\mathbb{N}_{-0}^n$ on $\hat{\mathcal{F}}$ and $S$ satisfies each $\tau_{ij}$-compatibility condition. Then $S$ is a passive collection
and the quotient algebra $\mathcal{F}/\langle O(S) \rangle$ is regular.

Proof. Since the partition $\{\Phi_\gamma\}_{\gamma \in \Gamma}$ of $\hat{\mathcal{F}}$ is compatible with the action of $\mathbb{N}_{-0}^n$ and every $f\in S$ is of the form (\ref{eq:f=u^i_alpha}), we have
$$ [D^\beta h^i_\alpha] \prec [u^i_{\alpha+\beta}] , \qquad \forall \beta\in \mathbb{N}_{-0}^n .$$
Any subset of $O(S)$ is therefore a conditionally solvable subset of $\mathcal{F}$.

Denote by $\Gamma_\gamma$ the set
$$\Gamma_\gamma = \{\gamma^{'}\in \Gamma: \theta \leq \gamma^{'} \leq \gamma\}.$$
It is easy to see that $C_{\gamma_1} \cap C_{\gamma_2} = \emptyset , \forall \gamma_1\neq \gamma_2 ,$
where $C_{\gamma_1}$, $C_{\gamma_2}$ given by (\ref{eq:Ogamma}). Thus the family $\{C_{\gamma} \}_{\gamma \in \Gamma_{\gamma_{*}}}$
 forms a partition of $O_{\gamma_{*}}$ (for all $\gamma_{*} \geq \theta$)
and 
\begin{equation} \label{eq:O_gamma}
O_{\gamma_{\star}} = C_{\gamma_{\star}} \cup   \bigcup_{\theta \leq \gamma < \gamma_{\star}} C_\gamma  .
\end{equation}

By the principle of trasfinite induction,  we shall prove that the ideal $\langle O_{\gamma} \rangle_{\mathcal{F}_{\gamma}}$ has a normalized set of generators $B_{\gamma}$ (for all $\gamma \geq \theta$) such that $B_{\gamma_1} \subset B_{\gamma_2}$ whenever $\gamma_1 < \gamma_2$.
 Assume $\gamma = \theta$, it follows that
$$ O_\theta = S\cap \mathcal{F}_\theta = S \cap \Phi_\theta = C_\theta . $$
Note that $C_\theta$ is a conditionally solvable subset of $\mathcal{F}$. We now distinguish two cases. In the first case, the leading terms of power series of $C_\theta$ are different.  Hence $C_\theta$ is a normalized set of generators of the ideal
$\langle O_\theta  \rangle_{\mathcal{F}_\theta} = \langle C_\theta  \rangle_{\mathcal{F}_\theta}$. In the second case,
some leading terms of power series of $C_\theta$ coincide. Let $f_i, f_j$ be two series in $C_\theta$ such that $lt f_i = lt f_j$.
From this it follows that $e_i - e_j$ is one of generators of $Syz\ lt S$ and either $f_i = f_j$ or
$f_i \neq f_j \wedge f_i - f_j \in \mathcal{F}_\gamma$ with $\gamma < \theta$. The last option is impossible because the pair of series $f_i, f_j$ must satisfy the compatibility condition (\ref{eq:compatib}), but $O_\gamma = \emptyset$ whenever $\gamma < \theta$.

Assume that for each $\gamma$ satisfying $\gamma_{*} > \gamma \geq \theta$, the ideal $\langle O_{\gamma} \rangle_{\mathcal{F}_{\gamma}}$ has a normalized set of generators $B_{\gamma}$.
It is now necessary to prove that there exists a normalized set of generators $B_{\gamma_{*}}$ of the ideal $\langle O_{\gamma_{*}} \rangle_{\mathcal{F}_{\gamma_{*}}}$.
We shall first prove that there are generators $Gen_{\gamma_{*}}$ of the ideal $\langle O_{\gamma_{*}} \rangle_{\mathcal{F}_{\gamma_{*}}}$ such that $Gen_{\gamma_{*}} \subseteq O_{\gamma_{*}}$ and the leading terms of power series of $Gen_{\gamma_{*}}$ are all distinct.

Because $O_{\gamma_{*}}$ may be represented by (\ref{eq:O_gamma}), it is useful to consider the set $C_{\gamma_{*}}$.
We study separately the two cases. In the first case, the leading terms of power series of $C_{\gamma_{*}}$ are different. Then
$$Gen_{\gamma_{*}} = (\bigcup_{\theta \leq \gamma < \gamma_{\star}} B_\gamma) \cup C_{\gamma_{\star}} .$$

We now suppose that some leading terms of power series of $C_{\gamma_{*}}$ coincide. Consider two series $f, g\in C_{\gamma_\star}$ such that $lt f = lt g$.  By definition of the orbit $O(S)$, there are $f_r, f_q \in S$ and $\mu, \nu \in \mathbb{N}^n$ satisfying
$$ f = D^\mu f_r, \qquad   g =  D^\eta g_q , $$
$$lt f = lt D^\mu f_r = D^\mu lt f_r = lt g = lt D^\eta g_q = D^\eta lt  g_q .  $$
The difference
$$d = X^\mu e_r - X^\eta e_q$$
is therefore a syzygy of $lt S$. It is not difficult to see that there exists $\sigma \in \mathbb{N}^n$ such that
\begin{equation}
\label{eq:d=} d = X^\sigma \tau_{rq} ,
\end{equation}
where $\tau_{rq}$ is one of the generators (\ref{eq:Syzygy}) of the module $Syz\ lt S$. Indeed, there are
$\sigma_1,\sigma_2,\sigma_3 \in \mathbb{N}^n$ satisfying
$$d - X^{\sigma_1} \tau_{rq} = (X^{\sigma_2} - X^{\sigma_3}) e_r .$$
The right side is a syzygy of $lt S$ if and only if $\sigma_2 = \sigma_3$. Thus the relation (\ref{eq:d=}) is true.

 Since $\tau_{rq}$ is one of the generators of the module $Syz\ lt S$ then, by assumption of our theorem, $S$ satisfies the $\tau_{rq}$-compatibility condition. This means
 $$D^{\alpha\diamond\beta}f_p - D^{\beta\diamond \alpha}f_q \in \Phi_{\gamma_1} \cap \langle O_{\gamma_1} \rangle_{\mathcal{F}_{\gamma_1}}$$
for some $\gamma_1 < \gamma_{*}$.
 Then applying $D^{\sigma}$ yields
\begin{equation} \label{eq:D^}
 f-g = D^{\alpha\diamond\beta +\sigma}f_p - D^{\beta\diamond \alpha+\sigma}f_q \in \Phi_{\gamma_2} \cap \langle O_{\gamma_2} \rangle_{\mathcal{F}_{\gamma_2}} , \qquad \gamma_1 < \gamma_2 < \gamma_{*},
\end{equation}
because the partition $\{\Phi_\gamma\}_{\gamma \in \Gamma}$ is compatible with the action of $\mathbb{N}_{-0}^n$.
We can conclude that $S$ satisfies the $d$-compatibility condition.

From (\ref{eq:D^}) we have
$$ f = g + \sum_{s=1,\dots,k, \, \nu\in \Delta } b^s_\nu D^\nu f_s ,$$
where $b^s_\nu,  D^\nu f_s \in \mathcal{F}_{\gamma_2}$, and $\Delta$ is a finite subset of $\mathbb{N}^n$. One may therefore include $g$ in
$Gen_{\gamma_{*}}$. In a similar way, one analyzes all pair $f, g \in C_{\gamma_{*}}$ such that $lt f = lt g$ and constructs the set $Gen_{\gamma_{*}}$ of generators of the ideal $\langle O_{\gamma_{*}} \rangle_{\mathcal{F}_{\gamma_{*}}}$.

Finally we can prove the existence of the normalized set of generators $B_{\gamma_{*}}$ of the ideal $\langle O_{\gamma_{*}} \rangle_{\mathcal{F}_{\gamma_{*}}}$. Let $f$ be any power series lying in $Gen_{\gamma_{*}} \cap \Phi_{\gamma_{*}}$. Then $f$ is written as
(\ref{eq:f=u^i_alpha}) with $h^i_\alpha \in \mathcal{F}_\gamma$, $\gamma < \gamma_{*}$. Suppose  $h^i_\alpha$ depends on
$v_1,\dots,v_l \in O(lt S) \cap \mathcal{F}_\gamma$. It follows from the induction assumption there are
$\tilde f_1,\dots,\tilde f_l \in B_\gamma$
such that $lt \tilde f_i = v_i$, $i=1,\dots,l$. Applying Lemma 1 to $h^i_\alpha$ we obtain
\begin{equation}
\label{eq:h=a1f1} h^i_\alpha = a_1 \tilde f_1 + \cdots + a_s \tilde f_s + r ,
\end{equation}
where $a_i \in \mathcal{F}_\gamma$ and $r$ does not depend on $v_1,\dots,v_l$. Then $\tilde f = u^i_\alpha + r$ is one of the generators of the ideal. In a similar way we can construct the normalized set of generators $B_{\gamma_{*}}$.

It is not difficult to show that
\begin{equation}
\label{eq:cupO} \bigcup_{\gamma \geq \theta} \langle  O_\gamma \rangle_{\mathcal{F}_\gamma} = \langle  O(S) \rangle .
\end{equation}
Indeed, any element $t\in \langle  O(S) \rangle$ can be written as a finite sum
$$ t = \sum_{i=1}^{r} a_i g_i    $$
with $a_i\in \mathcal{F}, g_i \in O(S)$. Observe that
$$ \mathcal{F} = \bigcup_{\gamma \in \Gamma} \mathcal{F}_\gamma , \qquad  O(S) = \bigcup_{\gamma \geq \theta} O_\gamma  $$
and the inclusions $ \mathcal{F}_\gamma \subset \mathcal{F}_{\gamma_{'}},  O_\gamma \subset O_{\gamma_{'}} $ hold for all
 $\gamma < \gamma_{'}$. Thus there is $\gamma_t\in \Gamma$ such that $a_i \in \mathcal{F}_{\gamma_t}$ and $g_i \in O_{\gamma_t}$, $i=1,\dots r$. The elements  $a_ig_i$ and  $t$  therefore belong to $\langle  O_{\gamma_t}  \rangle_{\mathcal{F}_{\gamma_t}}$.
It follows from  (\ref{eq:cupO}) that
\begin{equation}
\label{eq:Gen}   B = \bigcup_{\gamma \geq \theta} B_\gamma
\end{equation}
is the normalized set of generators of the ideal $\langle O(S) \rangle$. Then the quotient algebra $\mathcal{F}/\langle O(S) \rangle$ is regular, according to Proposition 1. The equality (\ref{eq:L=O}) holds by constructing the normalized set of generators.

\section{Conclusion}


This consideration leads to some important problems. If a collection $S$ is not passive then a natural question arises whenever there is a normalized set of generators of the ideal $\langle O(S) \rangle$. Janet \cite{Janet} described some method of constructing a passive system from given system of partial differential equations. Unfortunately, in general this method does not work.
 The second question is: How to choose a partition $\{\Phi_\gamma\}_{\gamma \in \Gamma}$ of $\hat{\mathcal{F}}$ to simplify  search of a passive collection for given differential ideal $\langle O(S) \rangle$? A similar problem is discussed in theory of Gr\"{o}bner bases \cite{Cox}.

\subsubsection*{Acknowledgments.} The research is supported in part by RFBR grants 13-01-00246, 12-01-00648 and by the Council on Grants of the President of the Russian Federation for Support of Leading Scientific Schools (Grants NSh-544.2012.1).

\end{document}